\documentclass[11pt,leqno]{amsart}
\usepackage{amsmath,amssymb,latexsym}
\usepackage{amscd}
\usepackage[OT4]{fontenc}
\usepackage{lmodern}
\usepackage{mathrsfs}
\usepackage{amsthm}
\usepackage{amsfonts}
\usepackage{url}
\usepackage[fleqn,tbtags]{mathtools}
\usepackage{enumerate}
\usepackage{graphicx}
\usepackage{color}
\usepackage[all]{xy}
\newcommand{\R}{\mathbb R}
\newcommand{\N}{\mathbb N}

\newcommand{\C}{\mathbb C}

\newcommand{\cB}{\mathcal B}

\newcommand{\cD}{\mathcal D}

\newcommand{\cK}{\mathcal K}
\newcommand{\cL}{\mathcal L}

\newcommand{\cS}{\mathcal S}
\newcommand{\im}{\operatorname{im}}

\newcommand{\ldss}{\mathcal{L}(s',s)}

\newcommand{\bproof}{{\raggedright\textbf{Proof.}} \ }
\newcommand{\proofmain}{{\raggedright\textbf{Proof of Theorem \ref{th_main}.}} \ }
\newcommand{\bqed}{\hspace*{\fill} $\Box $\medskip}

\let\epsilon\varepsilon
\let\phi\varphi
\let\rho\varrho

\theoremstyle{plain}
\newtheorem{Th}{Theorem}
\newtheorem{Cor}[Th]{Corollary}
\newtheorem{Lem}[Th]{Lemma}

\theoremstyle{definition}
\newtheorem{Def}[Th]{Definition}

\newtheorem*{Acknow*}{Acknowledgements}

\theoremstyle{remark}
\newtheorem{Rem}[Th]{Remark}

\title[Characterization of commutative algebras]{Characterization of commutative algebras embedded into the algebra of smooth operators}
\author{Tomasz Cia\'s}
\date{}

\begin{document}

\begin{abstract}
The paper deal with the noncommutative Fr\'echet ${}^*$-algebra $\ldss$ of the so-called smooth opertors, i.e. linear and continuous operators acting from the space $s'$ of slowly increasing 
sequences to the Fr\'echet space $s$ of rapidly decreasing sequences.  
By a canonical identification, this algebra of smooth operators can be also seen as the algebra of the rapidly decreasing matrices. 
We give a full description of closed commutative ${}^*$-subalgebras of this algebra and 
we show that every closed subspace of $s$ with basis is isomorphic (as a Fr\'echet space) to some closed commutative ${}^*$-subalgebra of $\ldss$. As a consequence, we give 
some equivalent formulation of the long-standing Quasi-equivalence Conjecture for closed subspaces of $s$.
\end{abstract}

\maketitle

\footnotetext[1]{


{The research of the author was supported by the National Center of Science, grant no. 2013/09/N/ST1/04410}.}

%\section{Introduction}

%We assume that $\lambda^\infty(A)\subset\ell_\infty$ to ensure that $\lambda^\infty(A)$ is an algebra.  

\section{Introduction}
In this paper we consider some specific noncommutative Fr\'echet algebra with involution -- known, for istance, as the algebra $\ldss$ of \emph{smooth operators} or as the algebra 
\[\cK_\infty:=\{(x_{j,k})_{j,k\in\N}\in\C^{\N^2}\colon\sup_{j,k\in\N}|x_{j,k}|j^qk^q<\infty\quad\textrm{for all }q\in\N_0\}\]
of \emph{rapidly decreasing matrices}.
Our main goal is to solve the following problems:
\begin{enumerate}
\item[(A)] characterize Fr\'echet ${}^*$-algebras isomorphic to closed commutative ${}^*$-subalgebras of $\ldss$ (Theorem \ref{th_main});
\item[(B)] characterize the underlying F\'echet spaces for the class of closed commutative ${}^*$-subalgebras of $\ldss$ (Theorem \ref{th_kothe_spaces_subalg_ldss}).
\end{enumerate}

By definition, $\ldss$ is the Fr\'echet ${}^*$-algebra consisting of all linear and continuous operators from the space
\begin{equation}\label{eq_s}
s':=\bigg\{(\xi_j)_{j\in\N}\in\C^\N\colon\sup_{j\in\N}|\xi_j|j^{-q}<\infty\quad\textrm{for some }q\in\N_0\} 
\end{equation}
of \emph{slowly increasing sequences} (this an inductive limit of Banach spaces) to the Fr\'echet space 
\begin{equation}\label{eq_s'}
s:=\{(\xi_j)_{j\in\N}\in\C^\N\colon\sup_{j\in\N}|\xi_j|j^q<\infty\quad\textrm{for all }q\in\N_0\} 
\end{equation}
of \emph{rapidly decreasing sequences}. 
Since there is a natural inclusion $\ldss$ into $\cB(\ell_2)$, the space $\ldss$ is indeed a ${}^*$-algebra with composition of operators and hilbertain involution taken from $\cB(\ell_2)$.
Futhermore, $\ldss$ with its natural topology of the uniform convergence on bounded subsets of $s'$ is a Fr\'echet space.

We prove, in particular, that the class of closed commutative ${}^*$-subalgebras of $\ldss$ is -- in the sense of an isomorphism of Fr\'echet ${}^*$-algebras -- 
exactly the class of K\"othe sequence algebras which are closed under taking square roots (Definition \ref{def_sqrt}) and at the same time isomorphic (as Fr\'echet spaces) to closed subspaces of $s$ (Theorem \ref{th_main}). 
This can be easily expressed by conditions on the corresponding K\"othe matrices. Moreover, from the solution of Problem (A) we get the solution of Problem (B): 
the underlying Fr\'echet spaces for closed commutative ${}^*$-subalgebras of $\ldss$ are just closed subspaces of $s$ with basis (Theorem \ref{th_kothe_spaces_subalg_ldss}). 

It has been already shown in \cite[Th. 4.8]{Cias} that every infinite-dimensional closed commutative ${}^*$-subalgebra $X$ of $\ldss$ is isomorphic as a Fr\'echet ${}^*$-algebra to a K\"othe sequence algebra
\begin{equation}\label{eq_intro_kothe_alg}
\lambda^\infty(||P_j||_q):=\bigg\{(\xi_j)_{j\in\N}\in\C^\N:\sup_{j\in\N}|\xi_j|\;||P_j||_q<\infty\text{ for all }q\in\N_0\bigg\},
\end{equation}
where $(P_j)_{j\in\N}$ is a suitable sequence of pairwise orthogonal, nonzero, self-adjoint projections belonging to $X$ and $(||\cdot||_q)_{q\in\N_0}$ is a sequence of continuous norms on $\ldss$ 
which determines the Fr\'echet space topology of $\ldss$. Hence, in order to solve problems (A) and (B), we needed to ``understand" matrices $(||P_j||_q)_{j\in\N,q\in\N_0}$,
where $(P_j)_{j\in\N}$ are pairwise orthogonal, nonzero, self-adjoint projections belonging to $\ldss$. 
The key step in solutions of our problems was to show that for every nuclear Fr\'echet space with the so-called dominating Hilbert norm $||\cdot||$ there is a topological embedding 
$V\colon X\hookrightarrow s$ such that $||Vx||_{\ell_2}=||x||$ for every $x\in X$ (Theorem \ref{th_top_embedding_DN}). 

In this context, it is worth mentioning that, by \cite[Th. 6.2]{Cias2}, a K\"othe sequence algebra of the form (\ref{eq_intro_kothe_alg}) is isomorphic to some closed ${}^*$-subalgebra of the algebra $s$ if 
and only if it is isomorphic as a Fr\'echet space to a complemented subspace of $s$. One would think that every closed commutative ${}^*$-subalgebra of $\ldss$ is of this type. 
Theorem \cite[Th. 6.9]{Cias2} shows that this is not true (in the proof we give a concrete contrexample) and the present paper shows that indeed -- taking into account the representation (\ref{eq_intro_kothe_alg}) --
the class of the underlying Fr\'echet spaces for closed commutative ${}^*$-subalgebras of $\ldss$ is as big as possible, i.e. as we mentioned above, it consists of all closed subspaces of $s$ with basis.

The algebra of smooth operators is a quite natural object. It can be seen as a kind of a noncommutative analogon of the commutative Fr\'echet algebra $s$ (look at the definition of $\cK_\infty$ above) and, 
in particular, it is isomorphic as a Fr\'echet space to $s$ (see \cite[Lemma 31.1]{MeV}) -- one of the most significant Fr\'echet spaces. 
On the level of Fr\'echet spaces, the space $s$ is isomorphic to many classical spaces of smooth functions
e.g. to the Schwartz space $\cS(\R)$ of smooth rapidly decreasing functions on the real line (which is a natural domain for the Fourier transform and differential operators with smooth coeficients), 
the space $C^\infty[-1,1]$ of smooth functions on the interval $[-1,1]$, the space $\cD(K)$ of smooth functions with support contained in a compact set $K\subset\R^n$, 
$\operatorname{int}(K)\neq\emptyset$, or the space $C^\infty(M)$ of smooth functions on an arbitrary compact smooth manifold $M$. 
Finally, the space $s$ carries all the information about nuclear locally convex spaces. Indeed, the K\=omura-K\=omura theorem states that a locally convex space 
is nuclear if and only if it is isomorphic to some closed subspace of a suitable big cartesian product of the space $s$ (see \cite[Th. 29.8]{MeV}). Consequently, a Fr\'echet space is nuclear 
(i.e. a Fr\'echet space in which every unconditionally convergent series is already absolutely convergent) if and only if it is isomorphic to some closed subspace of $s^\N$ (see \cite[Cor. 29.9]{MeV}).
 
Replacing $s'$ with the space $\cS'(\R)$ of tempered distributions and $s$ with $\cS(\R)$, we still end up with the same Fr\'echet ${}^*$-algebra -- this time consisting of
all linear and continuous operators acting from $\cS'(\R)$ to $\cS(\R)$, where ${}^*$-algebra operations are inherited from $L^2(\R)$. In fact, for every Fr\'echet space $E$ isomorphic to $s$ and its topological 
dual $E'$, the space $\cL(E',E)$ (with appropriate ${}^*$-algebra operations) of linear and continuous operators from $E'$ to $E$ is a Fr\'echet ${}^*$-algebra isomorphic to $\ldss$ (for details see 
\cite[Th. 1.10 \& Ex. 1.13]{Cias_phd} and \cite[Th. 2.1]{Dom}).  

The algebra $\ldss$ appears and plays a significant role in $K$-theory of Fr\'echet algebras 
(see Bhatt, Inoue and Ogi \cite[Ex. 2.12]{BhIO}, Cuntz \cite[p. 144]{Cuntz}, \cite[p. 64--65]{Cuntz2}, Gl\"ockner and Langkamp \cite{GlockLkamp}, 
Phillips \cite[Def. 2.1]{Phill}) and in $C^*$-dynamical systems (Elliot, Natsume and Nest \cite[Ex. 2.6]{ElNatNest}). Recently, Piszczek obtained several rusults concerning
closed ideals, automatic continuity, amenability, Jordan decomposition and Grothendieck-type inequality in $\cK_\infty$ 
(see Piszczek \cite{Pisz1,Pisz2,Pisz3,Pisz4} and his forthcoming paper ``The noncommutative Schwartz space is weakly amenable"). 

Recall that two bases $(x_j)_{j\in\N}$ and $(y_j)_{j\in\N}$ of a Fr\'echet space $X$ are called \emph{quasi-eqivalent} if  
there is a permutation $\sigma\colon\N\to\N$, a sequence $(\lambda_j)_{j\in\N}$ of nonzero scalars and a Fr\'echet space isomorphism $T\colon X\to X$ such that $Tx_j=\lambda_jy_{\sigma(j)}$ for $j\in\N$.
It is still an open question -- the so-called Quasi-equivalence Conjecture -- whether all bases in a nuclear Fr\'echet space are quasi-equivalent.
This conjecture was posed in 1961 by Mityagin \cite[\S 8.3]{Mit}. There were many attempts to solve this problem (see \cite{Zobin} and references therein), e.g. Crone and Robinson \cite{CR} showed that 
in every nuclear Fr\'echet space with the so-called regular basis all bases are quasi-equivalent. In particular, all bases in $s$ are quasi-equivalent but
the conjecture remains open for closed subspaces of $s$.
As a consequence of Theorem \ref{th_main}, we obtain some characterization of closed subspaces of $s$ which have all bases quasi-equivalent (Theorem \ref{th_quasi-eq}).
We hope that our characterization will shed some light on this long-standing conjecture.

%In order to solve these problems, we need to "understand" K\"othe matrices $(||P_k||_q)_{k\in\N,q\in\N_0}$, where 
%$(P_k)_{k\in\N}$ are pairwise orthogonal self-adjoint projections belonging to $\cL(s',s)$.

%\begin{enumerate}
% \item Characterize K\"othe matrices $A$ for which there is $(P_k)_{k\in\N}$ such that $\lambda^2(A)\cong\lambda^2(||P_k||_q)$ as Fr\'echet spaces.   
% \item Characterize K\"othe matrices $A$ for which there is $(P_k)_{k\in\N}$ such that 
% $\lambda^2(A)\cong\lambda^2(||P_k||_q)$ as Fr\'echet ${}^*$-algebras.
% $\lambda^2(A)=\lambda^2(||P_k||_q)$ as Fr\'echet ${}^*$-algebras.
% \begin{itemize}
% \item[($\alpha$)] $\forall q\in\N_0\;\exists r\in\N_0\;\exists C>0\;\forall k\in\N\quad a_{j,q}\leq C||P_k||_r$,
% \item[($\beta$)] $\forall r'\in\N_0\;\exists q'\in\N_0\;\exists C'>0\;\forall k\in\N\quad ||P_k||_{r'}\leq C'a_{j,q'}$.
%\end{itemize}
%\end{enumerate}

\section {Preliminaries}

In what follows, $\N$ will denote the set of natural numbers $\{1,2,\ldots\}$ and $\N_0:=\N\cup\{0\}$.

By $e_j$ we denote the vector in $\C^\N$ whose $j$-th coordinate equals 1 and the others equal 0.

By a \emph{Fr\'echet space} we mean a complete metrizable locally convex space over $\C$ (we will not use locally
convex spaces over $\R$). A \emph{Fr\'echet algebra} is a Fr\'echet space which is an algebra with continuous multiplication. 
A \emph{Fr\'echet ${}^*$-algebra} is a Fr\'echet algebra with continuous involution.

We use the standard notation and terminology. All the notions from functional analysis are explained in \cite{MeV}.
%and those from topological algebras in \cite{Fra} or \cite{Zel}.

We define the \emph{space of rapidly decreasing sequences} as the Fr\'echet space
\[s:=\bigg\{\xi=(\xi_j)_{j\in\N}\in\C^\N:|\xi|_q:=\bigg(\sum_{j=1}^\infty|\xi_j|^2j^{2q}\bigg)^{1/2}<\infty\text{ for all }q\in\N_0\bigg\}\]
with the topology corresponding to the system $(|\cdot|_q)_{q\in\N_0}$ of norms.
The strong dual of $s$ (i.e. the space of all continuous linear functionals on $s$ with the topology of uniform convergence on bounded
subsets of $s$, see e.g. \cite[Definition on p. 267]{MeV}) can be identified with the \emph{space of slowly increasing sequences} 
\[s':=\bigg\{\xi=(\xi_j)_{j\in\N}\in\C^\N:|\xi|_q':=\bigg(\sum_{j=1}^\infty|\xi_j|^2j^{-2q}\bigg)^{1/2}<\infty\text{ for some }q\in\N_0\bigg\}.\]
It is an easy exercise to show that replacing in the above definition of $s$ and $s'$ the $\ell_2$-norms by the corresponding $\sup$-norms -- as in (\ref{eq_s}) and (\ref{eq_s'}) -- or $\ell_p$-norms ($1\leq p<\infty$), 
we end up with the same sequence spaces.     

Closed subspaces of the space $s$ can be characterized by the so-called property (DN) (see \cite[Prop. 31.5]{MeV}). 
\begin{Def}\label{def_DN}
A Fr\'echet space $X$ with a fundamental system $(||\cdot||_q)_{q\in\N_0})$ of seminorms has the 
\emph{property} (DN)\index{property!(DN)} (see \cite[Def. on p. 359]{MeV}) if there is a continuous norm $||\cdot||$ on $X$ such that 
for all $q\in\N_0$ there is $r\in\N_0$ and $C>0$ such that 
\[||x||_q^2\leq C||x||\;||x||_r\]
for all $x\in X$. The norm $||\cdot||$ is called a \emph{dominating norm}\index{dominating norm}.

We will also use the condition (DN) in the following equivalent form (see \cite[Lemma 29.10]{MeV}):
there is a continuous norm $||\cdot||$ on $X$ such that for any $q\in\N_0$ and $\theta\in(0,1)$ there is $r\in\N_0$ and $C>0$ such that 
\[||x||_q\leq C||x||^{1-\theta}||x||_r^{\theta}\]
for all $x\in X$.
\end{Def}

The \emph{algebra of smooth operators} is defined as the Fr\'echet ${}^*$-algebra $\ldss$ of continuous linear operators from $s'$ to $s$ with the topology of uniform convergence on bounded sets in $s'$. 
The ${}^*$-algebra operations on $\ldss$ are inherited from $\cB(\ell_2)$ (note that every smooth operator is bounded on $\ell_2$).
The algebra $\ldss$ can be also seen as the algebra of matrices which are sometimes more handy; it is isomorphic, via the map $x\mapsto(\langle xe_k,ej\rangle)_{j,k\in\N}$, 
to the Fr\'echet ${}^*$-algebra
\[\cK_\infty:=\{(x_{j,k})_{j,k\in\N}\in\C^{\N^2}:\sup_{j,k\in\N}|x_{j,k}|j^qk^q<\infty\textrm{ for all }q\in\N_0\}\] 
of \emph{rapidly decreasing matrices} (with matrix multiplication and matrix complex conjugation). Therefore -- when reading this paper -- on can replace each ``$\ldss$" with ``$\cK_\infty$".
More details concerning $\ldss$ can be find in the introductions of \cite{Cias,Cias_phd,Cias2}.

\begin{Def}\label{def_Kothe}
A matrix $A=(a_{j,q})_{j\in\N,q\in\N_0}$ of nonnegative numbers such that
\begin{enumerate}
 \item[(i)] for each $j\in\N$ there is $q\in\N_0$ such that $a_{j,q}>0$
 \item[(ii)] $a_{j,q}\leq a_{j,q+1}$ for $j\in\N$ and $q\in\N_0$
\end{enumerate}
is called a \emph{K\"othe matrix}. For $1\leq p<\infty$ and a K\"othe matrix $A$ we define the \emph{K\"othe space} 
\[\lambda^p(A):=\bigg\{\xi=(\xi_j)_{j\in\N}\in\C^\N:|\xi|_{\lambda^p(A),q}^p:=\sum_{j=1}^\infty|\xi_j|^pa_{j,q}^p<\infty\text{ for all }q\in\N_0\bigg\}\]
and for $p=\infty$
\[\lambda^\infty(A):=\bigg\{\xi=(\xi_j)_{j\in\N}\in\C^\N:|\xi|_{\lambda^\infty(A),q}:=\sup_{j\in\N}|\xi_j|a_{j,q}<\infty\text{ for all }q\in\N_0\bigg\}\]
with the locally convex topology given by the seminorms $(|\cdot|_{\lambda^p(A),q})_{q\in\N_0}$ (see e.g. \cite[Definition p. 326]{MeV}). 
For simplicity, we will also write $\lambda^p(a_{j,q})$ (i.e. only the entries of the matix) instead of $\lambda^p(A)$.
\end{Def}

All spaces $\lambda^p(A)$ are Fr\'echet spaces (\cite[Lemma 27.1]{MeV}).
By the Dynin-Mityagin Theorem (see \cite[Th. 28.12]{MeV}) if $\lambda^p(A)$ is nuclear (i.e. for all $q$ there is $r$ such that $\sum_{j=1}^\infty\frac{a_{j,q}}{a_{j,r}}<\infty$), 
then the sequence of vectors $(e_j)_{j\in\N}$ is an absolute Schauder basis of $\lambda^p(A)$. 
We will also use the following result: $\lambda^p(A)$ is nuclear for some $1\leq p\leq\infty$ if and only if $\lambda^p(A)=\lambda^q(A)$ as Fr\'echet spaces for all $1\leq p,q\leq\infty$ (\cite[Prop. 28.16]{MeV}).

For K\"othe matrices $A=(a_{j,q})_{j\in\N,q\in\N_0}$ and $B=(b_{j,q})_{j\in\N,q\in\N_0}$ such that 
\[\forall q\in\N_0\;\exists r\in\N_0\;\exists C>0\;\forall j\in\N\quad a_{j,q}\leq Cb_{j,r},\]
we write $A\prec B$.
If $A\prec B$ and $B\prec A$, then we write $A\sim B$. For a bijection $\sigma\colon\N\to\N$, $A_\sigma$ denotes the K\"othe matrix $(a_{\sigma(j),q})_{j\in\N,q\in\N_0}$ and, moreover, $A^2$ is by definition 
the K\"othe matrix $(a_{j,q}^2)_{j\in\N,q\in\N_0}$.

By \cite[Prop. 3.1]{Pir}, $\lambda^1(A)$ with pointwise multiplication is an algebra (and, clearly, if $\lambda^1(A)$ is nuclear so is $\lambda^p(A)$ for $1<p\leq\infty$) if and only if $A\prec A^2$; in this case $\lambda^p(A)$ is called 
a \emph{K\"othe algebra}. 

\begin{Def}\label{def_sqrt}
We say that a K\"othe space $\lambda^p(A)$ is \emph{closed under taking square root} if it is a K\"othe algebra and for each sequence $(\xi_j)_{j\in\N}$ in $\lambda^p(A)$ of nonnegative scalars, 
the sequence $(\sqrt{\xi_j})_{j\in\N}$ belongs to $\lambda^p(A)$ as well.
\end{Def}

\section{Commutative subalgebras of $\ldss$}

Our main result -- solving Problem (A) from Introduction -- reads as follows:

\begin{Th}\label{th_main}
For a K\"othe matrix $A=(a_{j,q})_{j\in\N,q\in\N_0}$ the following assertions are equivalent: 
\begin{enumerate}
 \item[\upshape(1)] $\lambda^2(A)$ is isomorphic as a Fr\'echet ${}^*$-algebra to some closed commutative ${}^*$-subalgebra of $\ldss$;
 \item[\upshape(2)] there is a sequence $(P_j)_{j\in\N}$ of nonzero pairwise orthogonal self-adjoint projections belonging to $\cL(s',s)$ such that 
 $A\sim(||P_j||_q)_{j\in\N,q\in\N_0}$;
 \item[\upshape(3)] there is a basic sequence $(f_j)_{j\in\N}$ of $s$ which is at the same time an orthonormal sequence in $\ell_2$ such that 
 $A\sim(|f_j|_q)_{j\in\N,q\in\N_0}$;
 \item[\upshape(4)] the matrix $A$ satisfies the following conditions:  
 \begin{enumerate}
  \item[\upshape(i)] $\forall q\in\N_0\;\exists r\in\N_0\quad \sum_{j=1}^\infty\frac{a_{j,q}}{a_{j,r}}<\infty$,
  \item[\upshape(ii)] $\exists p\in\N_0\;\exists C>0\;\forall j\in\N\quad a_{j,p}\geq\frac{1}{C}$,
  \item[\upshape(iii)] $A^2\prec A$;
 \end{enumerate}
 \item[\upshape(5)] the matrix $A$ satisfies the following conditions: 
 \begin{enumerate}
  \item[\upshape(i)] $\forall q\in\N_0\;\exists r\in\N_0\quad \sum_{j=1}^\infty\frac{a_{j,q}}{a_{j,r}}<\infty$,
  \item[\upshape(ii)] $\exists p\in\N_0\;\forall q\in\N_0\;\exists r\in\N_0,C>0\;\forall j\in\N\quad a_{j,q}^2\leq Ca_{j,p}a_{j,r}$,
  \item[\upshape(iii)] $A\sim A^2$;
 \end{enumerate}
 \item[\upshape(6)] $\lambda^2(A)$ is a nuclear biprojective K\"othe algebra with the property \emph{(DN)};
 \item[\upshape(7)] $\lambda^2(A)$ is isomorphic as a Fr\'echet space to a closed subspace of $s$ and it is closed under taking square roots;
 \item[\upshape(8)] $\lambda^2(A)$ is a nuclear Fr\'echet space with $||\cdot||_{\ell_2}$ as a dominating norm.
\end{enumerate}
\end{Th}

\begin{Rem}
By \cite[Th. 4.8]{Cias}, every infinite-dimensional closed commutative ${}^*$-subalgebra of $\ldss$ is isomorphic as a Fr\'echet ${}^*$-algebra to some K\"othe sequence algebra, and thus the above theorem
characterizes in fact all infinite-dimensional closed commutative ${}^*$-subalgebras of $\ldss$.
\end{Rem}

\begin{Rem} 
Homology theory is out of the mainstream of this paper, and therefore we only recall that the a Fr\'echet algebra $A$ is called \emph{biprojective} if 
the product map $\pi\colon A\hat\otimes A\to A$, $a\otimes b\mapsto ab$, has a right inverse in the category of $A$-bimodules (see also \cite[Def. IV.5.1]{Hel}).
The characterization of biprojectivity for K\"othe algebras in terms of the corresponding K\"othe matrices is due to Pirkovskii \cite[Th. 5.9]{Pir}. 
\end{Rem}

The following theorem leads to the solution of Problem (B) (see Introduction).

\begin{Th}\label{th_kothe_spaces_subalg_ldss}
For a K\"othe matrix $A=(a_{j,q})_{j\in\N,q\in\N_0}$ the following assertions are equivalent:
\begin{enumerate}
 \item[\upshape(1)] $\lambda^2(A)$ is isomorphic as a Fr\'echet space to some closed commutative ${}^*$-subalgebra of $\ldss$; 
 \item[\upshape(2)] the matrix $A$ satisfies the following conditions: 
 \begin{enumerate}
  \item[\upshape(i)] $\forall q\in\N_0\;\exists r\in\N_0\quad \sum_{j=1}^\infty\frac{a_{j,q}}{a_{j,r}}<\infty$,
  \item[\upshape(ii)] $\exists p\in\N_0\;\forall q\in\N_0\;\exists r\in\N_0,C>0\;\forall j\in\N\quad a_{j,q}^2\leq Ca_{j,p}a_{j,r}$;
 \end{enumerate}
 \item[\upshape(3)] $\lambda^2(A)$ is a nuclear Fr\'echet space with the property $\mathrm{(DN)}$.
 \item[\upshape(4)] $\lambda^2(A)$ is a Fr\'echet space isomorphic to a closed subspace of $s$.
\end{enumerate}
\end{Th}
\bproof
It is well known that a K\"othe space $\lambda^2(A)$ is nuclear and has the property (DN) if and only if $A$ satifies condtions (i) and (ii),
respectively. This shows (1)$\Rightarrow$(2)$\Rightarrow$(3), because $\ldss\cong s$ as Fr\'echet spaces and nuclearity together 
with the property (DN) are inherited by closed subspaces.

In order to prove (3)$\Rightarrow$(1), choose $p_0\in\N_0$ for which $|\cdot|_{\lambda^2(A),p_0}$ is dominating norm on $\lambda^2(A)$ 
and define the K\"othe matrix $B:=(a_{j,q}/a_{j,p_0})_{j\in\N,q\in\N_0}$. 
Let $T\colon\lambda^2(B)\to \lambda^2(A)$, $T\xi:=(\xi_j/a_{j,p_0})_{j\in\N}$. Then, for $\xi\in\lambda^2(B)$ and $q\in\N_0$, we have
$|T\xi|_{\lambda^2(A),q}=|\xi|_{\lambda^2(B),q}$,
and thus $T$ is an isomorphism of Fr\'echet spaces. Moreover, for $q\in\N_0$ there is $r\in\N_0$ and $C>0$ such that 
\[|\xi|_{\lambda^2(B),q}^2=|T\xi|_{\lambda^2(A),q}^2\leq C |T\xi|_{\lambda^2(A),p_0}|T\xi|_{\lambda^2(A),r}
=C||\xi||_{\ell_2}|\xi|_{\lambda^2(B),r}\]
$\xi\in\lambda^2(B)$.
Therefore, $||\cdot||_{\ell_2}$ is a dominating norm on $\lambda^2(B)$, and the conclusion follows from Theorem \ref{th_main} 
(implication (8)$\Rightarrow$(1)). 

Finally, equivalence (3)$\Leftrightarrow$(4) follows directly from \cite[Prop. 31.5]{MeV}.
\bqed

As a direct consequence of Theorems \ref{th_main}, \ref{th_kothe_spaces_subalg_ldss} and 
\cite[Cor. 28.13, Prop. 28.16, Prop. 31.5]{MeV} we also get the following characterization.

\begin{Cor}\label{cor_frechet_spaces_subalg_ldss} 
For a Fr\'echet space $E$ the following assertions are equivalent:
\begin{enumerate}
 \item[\upshape(1)] $E$ is isomorphic as a Fr\'echet space to some closed commutative ${}^*$-subalgebra of $\cL(s',s)$;
 \item[\upshape(2)] $E$ is isomorphic as a Fr\'echet space to a K\"othe space $\lambda^2(|f_j|_q)$ 
 for some basic sequence $(f_j)_{j\in\N}$ of $s$ which is at the same time an orthonormal sequence in $\ell_2$;
 \item[\upshape(3)] $E$ is nuclear, has a basis and the property $\mathrm{(DN)}$;
 \item[\upshape(4)] $E$ is isomorphic as a Fr\'echet space to some closed subspace of $s$ with basis.
\end{enumerate}
\end{Cor}

Now, we shall provide the proof of Theorem \ref{th_main}. We start with some lemmas involving short exact sequences of Fr\'echet-Hilbert spaces.

\begin{Lem}\label{lem_EGF}
Let  
\[0\longrightarrow E\stackrel{j}\longrightarrow F\stackrel{q}\longrightarrow G\longrightarrow 0\]
be a short exact sequence of Fr\'echet-Hilbert spaces. Let $||\cdot||_E$ be a continuous Hilbert norm on $E$ and let $||\cdot||_G$ be a continuous norm on $G$.
Then there is a continuous Hilbert norm $||\cdot||_F$ on $F$ with 
\begin{enumerate}
 \item[\upshape($\alpha$)] $||x||_E\sim ||j(x)||_F$ for $x\in E$;
 \item[\upshape($\beta$)] $||q(x)||_G\leq\inf_{y\in \ker q}||x-y||_F$ for $x\in F$.
\end{enumerate}
\end{Lem}

\bproof
Without lost of generality we may assume that $E$ is a closed subspace of $F$ and $j(x)=x$ for $x\in E$.
%Since $j(E)=\ker q$ is a Fr\'echet space, the open mapping theorem implies that the map $j^{-1}\colon j(E)\to E$ is continuous. Hence, we can find a continuous seminorm $||\cdot||_{j(E)}$ on $j(E)$ such that $||x||_E\leq ||j(x)||_{j(E)}$ for $x\in E$. 
By \cite[Remark 22.8]{MeV}, there is a continuous Hilbert norm $||\cdot||_1$ on $F$ with $||x||_E\leq||x||_1$ for $x\in E$. Moreover, by the continuity of $q$, there is a continuous Hilbert norm 
$||\cdot||_2$ on $F$ such that 
\[||q(x)||_G\leq\inf_{z\in E}||x-z||_2\]
for $x\in F$.
Then $||\cdot||:=(||\cdot||_1^2+||\cdot||_2^2)^{1/2}$ is a continuous Hilbert norm on $F$ such that
\begin{enumerate}
 \item[\upshape(i)] $||x||_E\leq||x||$ for $x\in E$;
 \item[\upshape(ii)] $||q(x)||_G\leq\inf_{z\in E}||x-z||$ for $x\in F$.
\end{enumerate}
Define $||\cdot||_F$, by
\[||x||_F^2=\inf_{z\in E}\big(||z||_E^2+||x-z||^2\big)\]
for $x\in F$.

To show condition ($\alpha$), let us fix $x\in E$. Taking $z:=x$ in the infimum, we get
\[||x||_F\leq ||x||_E.\]
Moreover, for $z\in E$, we obtain
\[||x||_E^2\leq (||z||_E+||x-z||_E)^2=||z||_E^2+||x-z||_E^2+2||z||_E\,||x-z||_E.\]
If $||z||_E\leq||x-z||_E$, then $||z||_E\,||x-z||_E\leq||x-z||_E^2$. Otherwise, $||z||_E\,||x-z||_E\leq||z||_E^2$ so, by (i),
\[||x||_E^2\leq 3\big(||z||_E^2+||x-z||_E^2\big)\leq 3\big(||z||_E^2+||x-z||^2\big).\]
In consequence,
\begin{equation}\label{eq_E_leq_sqrt3F}
||x||_E\leq \sqrt{3}\inf_{z\in E}\big(||z||_E^2+||x-z||^2\big)^{1/2}=\sqrt{3}||x||_F,
\end{equation}
and thus condition ($\alpha$) if fulfilled.
By (ii), for $x\in F$, we get
\begin{align*}
||q(x)||_G^2&\leq\inf_{z\in E}||x-z||^2=\inf_{z_1\in E}\inf_{z_2\in E}\big(||z_2||_E^2+||x-z_1||^2\big)=\inf_{z_1,z_2\in E}\big(||z_2||_E^2+||x-(z_1+z_2)||^2\big)\\
&=\inf_{z_1\in E}\inf_{z_2\in E}\big(||z_2||_E^2+||(x-z_1)-z_2||^2\big)=\inf_{z\in E}||x-z||_F^2,
\end{align*}
which yields condition ($\beta$). 

Now define a Hilbert norm $||(\cdot,\cdot)||_H$ on $H:=E\times F$ by
\[||(x,y)||_H^2:=||x||_E^2+||y||^2.\]
Let 
\[R\colon (H,||(\cdot,\cdot)||_H)\to F,\quad (x,y)\mapsto x+y\]
and let
\[Q\colon (H,||\cdot||_H)\to(H/\ker R, ||\cdot||_{H/\ker R})\]
be the quotient map.
Then 
\begin{align}
\label{eq_H/kerR=F}
\begin{split}
||Q(x,y)||_{H/\ker R}&=\inf_{(x',y')\in \ker R}||(x,y)-(x',y')||_H=\inf_{x'\in E}||(x,y)-(x',-x')||_H\\
&=\inf_{x'\in E}\big(||x-x'||_E^2+||y+x'||^2\big)^{1/2}=\inf_{z\in E}\big(||z||_E^2+||(x+y)-z||^2\big)^{1/2}\\
&=||x+y||_F.
\end{split}
\end{align}
In particular
\begin{equation}\label{eq_F=H/kerR}
||\cdot||_F=||Q(0,\cdot)||_{H/\ker R}, 
\end{equation}
hence $||\cdot||_F$ is a seminorm. Moreover, since $||\cdot||_F\leq||\cdot||$, it follows that $||\cdot||_F$ is continuous.

Now assume that $||x||_F=0$ for some $x\in F$. Then, by condition ($\beta$), $||q(x)||_G\leq ||x||_F=0$, and since $||\cdot||_G$ is a norm, we have
$x\in\ker q=E$. Moreover, by (\ref{eq_E_leq_sqrt3F}), $||x||_E\leq\sqrt{3}||x||_F=0$, but $||\cdot||_E$ is a norm, so $x=0$. This shows that $||\cdot||_F$ is a norm, and therefore, 
by (\ref{eq_H/kerR=F}), $||\cdot||_{H/\ker R}$ is a norm as well (i.e. $\ker R$ is a closed subspace of $H$ -- see e.g. \cite[Lemma 5.10]{MeV}). 
Since the quotient norm induced by a Hilbert norm is a Hilbert norm, the norm $||\cdot||_{H/\ker R}$ is hilbertian (see e.g. \cite[p. 44]{Schr}). Therefore, by (\ref{eq_F=H/kerR}),
$||\cdot||_F$ satisfies the parallelogram law, hence $||\cdot||_F$ is a Hilbert norm.
\bqed

%\begin{Lem}
%Let $E$ be a closed subspace of $s$ with a dominating norm $||\cdot||_E$ and let $||\cdot||_s$ be a norm on $s$ such that $U_s\cap E=U_E$
%and $||\cdot||_{q(U_s)}$ is a dominating norm on $s/E$, where $U_s:=\{x\in s\colon ||x||_s<1\}$, $U_E:=\{x\in E\colon ||x||_E<1\}$ and $q\colon s\to s/E$ is the quotient map.
%Then $||\cdot||_s$ is a dominating norm on $s$.
%\end{Lem}

%\begin{Lem}\label{lem_||cdot||_s_DN}
%Let $E$ be a closed subspace of $s$ with a dominating norm $||\cdot||_E$ and let $||\cdot||_s$ be a continuous seminorm (??) on $s$ such that 
%${||\cdot||_s}_{\mid E}=||\cdot||_E$ and $||\cdot||_{s/E}$, $||q(x)||_{s/E}:=\inf_{y\in E}||x-y||_s$, is a dominating norm on $s/E$ (here $q\colon s\to s/E$ is the quotient map).
%Then $||\cdot||_s$ is a dominating norm on $s$.
%\end{Lem}
%\bproof
%Use \cite[Lemma 4.4]{V4}
%\bqed

An easy proof of our next lemma is left to the reader.
\begin{Lem}\label{lem_completionXY}
Let us consider the diagram 
\begin{displaymath}
\xymatrix{
(Y,||\cdot||_Y)\ar @{^{(}->}[r]^\iota \ar @{^{(}->}[d]^{j_Y}         & (X,||\cdot||_X) \ar @{^{(}->}[d]^{j_X}\\
(\widetilde{Y},||\cdot||_{\widetilde{Y}}) \ar[r]^{\widetilde{\iota}} & (\widetilde{X},||\cdot||_{\widetilde{X}})  
}
\end{displaymath}
where 
\begin{itemize}
\item $(X,||\cdot||_X)$ and $(Y,||\cdot||_Y)$ are normed spaces such that $Y$ is a closed subspace of $(X,||\cdot||_X)$ and $||\cdot||_Y\sim{||\cdot||_X}_{\mid Y}$;
\item $(\widetilde{Y},||\cdot||_{\widetilde{Y}})$ and $(\widetilde{X},||\cdot||_{\widetilde{X}})$ are the completions of $(Y,||\cdot||_Y)$ and $(X,||\cdot||_X)$, respectively;
\item $\iota$, $j_Y$ and $j_X$ are the canonical inclusions;
\item $\widetilde{\iota}$ is the continuous linear extensions of the map $j_X\circ\iota\colon Y\to \widetilde{X}$.
\end{itemize}
 
Then $\widetilde\iota(\widetilde{Y})$ is a closed subspace of $(\widetilde{X},||\cdot||_{\widetilde{X}})$.
\end{Lem}

%\bproof
%Let $\widehat{Y}:=\widetilde\iota(\widetilde{Y})$.
%We shall show that $\widehat{Y}=\overline{j_X(Y)}^{||\cdot||_{\widetilde{X}}}$. 
%For $z\in\widehat{Y}$ there is a sequence 
%$(y_n)_{n\in\N}\subset Y$ such that $(j_Y(y_n))_{n\in\N}$ is convergent in $(\widetilde{Y},||\cdot||_{\widetilde{Y}})$ and
%\[z=\widetilde{\iota}\big(\lim_{n\to\infty}j_Y(y_n)\big)=\lim_{n\to\infty}\widetilde{\iota}(j_Y(y_n))=\lim_{n\to\infty}j_X(y_n),\]
%where the second identity follows from the continuity of $\widetilde\iota$. 
%Therefore, $\widehat{Y}\subseteq\overline{j_X(Y)}^{||\cdot||_{\widetilde{X}}}$. 

%In order to prove the opposite inclusion, take $z\in\overline{j_X(Y)}^{||\cdot||_{\widetilde{X}}}$. 
%Then $z=\lim_{n\to\infty}j_X(y_n)$ for some sequence $(y_n)_{n\in\N}\subset Y$. Since $||\cdot||_Y\sim{||\cdot||_X}_{\mid Y}$, there is $C>0$ such that
%\[||y_n-y_m||_Y\leq C||y_n-y_m||_X=C||j_X(y_n)-j_X(y_m)||_{\widetilde{X}},\]
%and thus $(y_n)_{n\in\N}$ is a Cauchy sequence in $Y$. Consequently, $(j_Y(y_n))_{n\in\N}$ is convergent in $\widetilde{Y}$ and as above
%\[z=\lim_{n\to\infty}j_X(y_n)=\lim_{n\to\infty}\widetilde{\iota}(j_Y(y_n))=\widetilde{\iota}\big(\lim_{n\to\infty}j_Y(y_n)\big)\in\widehat{Y}.\]
%Hence,  $\overline{j_X(Y)}^{||\cdot||_{\widetilde{X}}}\subseteq\widehat{Y}$.
%\bqed

%Now we prove a key lemma which shows that -- under minor conitions -- one can extend a dominating Hilbert norm from a closed subspace of $s$ to a dominating Hilbert norm on the whole space $s$. 

\begin{Lem}\label{lem_extension_of_Hilbert_norm}
Let  
\[0\longrightarrow E\stackrel{j}\longrightarrow F\stackrel{q}\longrightarrow G\longrightarrow 0\]
be a short exact sequence of Fr\'echet-Hilbert spaces. Let $||\cdot||_E$ be a continuous Hilbert norm on $E$ and let us assume that $G$ has a continuous norm.
Then there is a continuous Hilbert norm $||\cdot||$ on $F$ such that $||j(x)||=||x||_E$ for $x\in E$. 
\end{Lem}
\bproof
Without lost of generality we may assume that $E$ is a closed subspace of $F$ and $j(x)=x$ for $x\in E$.
Let $|\cdot|_G$ be a continuous norm on $G$. By Lemma \ref{lem_EGF}, there is a continuous Hilbert norm $||\cdot||_F$ on $F$ such that
\begin{enumerate}
 \item[\upshape($\alpha$)] $||\cdot||_E\sim {||\cdot||_F}_{\mid{E}}$ for $x\in E$;
 \item[\upshape($\beta$)] $|q(x)|_G\leq\inf_{y\in \ker q}||x-y||_F=:||x||_G$ for $x\in F$.
\end{enumerate}
Consider the diagram
\begin{displaymath}
\xymatrix{
0\ar[r] & (E,||\cdot||_E)\ar @{^{(}->}[r]^j \ar @{^{(}->}[d]^{j_E}         & (F,||\cdot||_F) \ar@{^{(}->}[d]^{j_F}\ar[r]^q                                   & (G,||\cdot||_G)\ar[r] & 0\\
{}      & (\widetilde{E},||\cdot||_{\widetilde{E}}) \ar[r]^{\widetilde{j}} & (\widetilde{F},||\cdot||_{\widetilde{F}})\ar@<1ex>[dl]^{\pi_0}\ar@<1ex>[l]^\pi  & {}                    &  \\
{}      & (\widehat{E},|\cdot|_{\widetilde{F}})\ar@<1ex>[u]_{\kappa}       & {}                                                                              & {}                    &  \\
}  
\end{displaymath}
where
\begin{itemize}
 \item $(\widetilde{E},||\cdot||_{\widetilde{E}})$ and $(\widetilde{F},||\cdot||_{\widetilde{F}})$ are the completions of $(E,||\cdot||_E)$ and $(F,||\cdot||_F)$, respectively;
 \item $j_E$ and $j_F$ are the canonical inclusions;
 \item $\widetilde{j}$ is the continuous extensions of the map $j_F\circ j\colon E\to \widetilde{F}$. 
\end{itemize}
Since, by ($\beta$), $||\cdot||_{G}$ is a norm, $E$ is a closed subspace of $(F,||\cdot||_F)$, and thus, by Lemma \ref{lem_completionXY}, $\widehat{E}:=\widetilde{j}(\widetilde{E})$ is a closed subspace of 
the Hilbert space $(\widetilde{F},||\cdot||_{\widetilde{F}})$. Then $\pi_0$ is defined to be the orthogonal projections onto $\widehat{E}$, $\kappa$ is the inverse of $\widetilde\iota$ on $\widehat{E}$ 
and $\pi:=\kappa\circ\pi_0$.

Let us define
\[||x||^2:=||(\pi\circ j_F)(x)||_{\widetilde{E}}^2+||q(x)||_G^2\]
for $x\in F$. Clearly, $||\cdot||$ is a continuous seminorm on $F$. 
Next, for $x\in E$, we get
\[||x||=||(\pi\circ j_F)(x)||_{\widetilde{E}}=||j_E(x)||_{\widetilde{E}}=||x||_E,\]
hence ${||\cdot||}_{\mid{E}}=||\cdot||_E$. For $x\in F$ with $||x||=0$ we have $||q(x)||_G=0$, and thus $x\in E$. Then we have also
$||x||_E=||x||=0$ so $x=0$. This shows that $||\cdot||$ is a norm.

Finally, since $||\cdot||_{\widetilde{E}}$ and $||\cdot||_G$ are Hilbert norms, 
$||\cdot||$ is a Hilbert norm as well.
\bqed

The proof of the following lemma is a slight modification of the Vogt's proof of \cite[Lemma 3.1]{V6}. For the convenience of the reader we present here the full argument.    
 
\begin{Lem}\label{lem_extenstion_of_D-Hilbert_N}
Let  
\[0\longrightarrow E\stackrel{j}\longrightarrow F\stackrel{q}\longrightarrow G\longrightarrow 0\]
be a short exact sequence of Fr\'echet-Hilbert spaces. Let $||\cdot||_E$ be dominating Hilbert norm on $E$ and let us assume that $G$ has the property \emph{(DN)}.
Then there is a dominating Hilbert norm $||\cdot||$ on $F$ such that $||j(x)||=||x||_E$ for $x\in E$. 
\end{Lem}
\bproof
Without lost of generality we may assume that $E$ is a closed subspace of $F$ and $j(x)=x$ for $x\in E$.
By Lemma \ref{lem_extension_of_Hilbert_norm}, there is continuous Hilbert norm $||\cdot||_F$ on $F$ such that ${||\cdot||_F}_{\mid{E}}=||\cdot||_E$ and, by assumption,
there is a dominating Hilbert norm $||\cdot||_G$ on $G$.  Define
\[||x||^2:=||x||_F^2+||qx||^2_G\]
for $x\in F$ and let $U_0:=\{x\in E\colon ||x||<1\}$. Then, clearly, $||\cdot||$ is a continuous Hilbert norm on $F$. We shall show that $||\cdot||$ is a dominating norm.

For an absolutely convex zero neighborhood $U$ in $F$, let $||\cdot||_U$ denote the Minkowski functional for $U$ and let $||\cdot||_{\widehat{U}}\colon G\to[0,\infty)$ be defined by
$||qx||_{\widehat{U}}:=\inf\{||x+y||_U\colon y\in E\}$ for $x\in F$.

Since $||\cdot||_E$ is a dominating norm on $E$, for each absolutely convex zero neighborhood $U$ in $F$ there is an absolutely convex zero neighborhood $W$ in $F$ with $W\subseteq U\cap U_0$ such that
\begin{equation}\label{eq_DN_E}
||z||_U\leq||z||_E^{2/3}||z||_W^{1/3},
\end{equation}
for all $z\in E$. 
Since $||\cdot||_G$ is a dominating norm on $G$, there is an absolutely convex zero neighborhood $V$ in $F$ with $V\subseteq W$ such that
\begin{equation}\label{eq_DN_G}
||qx||_{\widehat{W}}\leq||qx||_G^{3/4}||qx||_{\widehat{V}}^{1/4}\leq||x||^{3/4}||x||_V^{1/4}. 
\end{equation}
for all $x\in F$.

Let us fix $x\in F$ and $\epsilon>0$. By the very definition of the norm $||\cdot||_{\widehat{W}}$, there is $y\in F$ such that 
\[qx=qy \quad\text{and}\quad ||y||_W\leq||qx||_{\widehat{W}}+\epsilon,\]
and thus (\ref{eq_DN_G}) yields
\begin{equation}\label{eq_yW_leqx+epsilon}
||y||_W\leq||x||^{3/4}||x||_V^{1/4}+\epsilon.
\end{equation}
Hence, for $z:=x-y$, we have 
\begin{align*}
||z||_E\leq||x||_F+||y||_F\leq||x||+||y||_W\leq||x||+||x||^{3/4}||x||_V^{1/4}+\epsilon,
\end{align*}
and since
\[||x||=||x||^{3/4}||x||^{1/4}\leq||x||^{3/4}||x||_V^{1/4},\] 
it follows that
\begin{equation}\label{eq_zE}
||z||_E\leq2||x||^{3/4}||x||_V^{1/4}+\epsilon.
\end{equation}
Moreover, by (\ref{eq_yW_leqx+epsilon}), 
\[||z||_W\leq||x||_W+||y||_W\leq||x||_V+||x||^{3/4}||x||_V^{1/4}+\epsilon\leq2||x||_V+\epsilon.\]
Now, since for $a,b\geq0$ we have $(a+b)^{2/3}\leq3(a^{2/3}+b^{2/3})$, it follows from (\ref{eq_DN_E}) and (\ref{eq_zE}) that
\begin{align}\label{eq_zU}
\begin{split}
||z||_U&\leq\left(2||x||^{3/4}||x||_V^{1/4}+\epsilon\right)^{2/3}\left(2||x||_V+\epsilon\right)^{1/3}\\
&\leq\left[3\left(2||x||^{3/4}||x||_V^{1/4}\right)^{2/3}+\epsilon^{2/3}\right]\left[\left(2||x||_V\right)^{1/3}+\epsilon^{1/3}\right]\\
&=\left(3\sqrt[3]{4}||x||^{1/2}||x||_V^{1/6}+3\epsilon^{2/3}\right)\left(\sqrt[3]{2}||x||_V^{1/3}+\epsilon^{1/3}\right)=6||x||^{1/2}||x||_V^{1/2}+h(\epsilon),
\end{split}
\end{align}
where $h(\epsilon)\searrow0$ when $\epsilon \searrow0$.

Next, by (\ref{eq_yW_leqx+epsilon}),
\begin{equation}
||y||_U\leq||y||_W\leq||x||^{3/4}||x||_V^{1/4}+\epsilon=||x||^{1/2}||x||_V^{1/2}\left(\frac{||x||}{||x||_V}\right)^{1/4}+\epsilon\leq||x||^{1/2}||x||_V^{1/2}+\epsilon.
\end{equation}
Consequently, by (\ref{eq_zU}),
\[||x||_U\leq||y||_U+||z||_U\leq 7||x||^{1/2}||x||_V^{1/2}+h(\epsilon)+\epsilon,\]
and finally, since $h(\epsilon)+\epsilon$ can be chosen arbitrary small, we have
\[||x||_U^2\leq 49||x||\;||x||_V,\]
which shows that $||\cdot||$ is a dominating norm on $F$.
\bqed

\begin{Cor}\label{cor_extension_hilbertDN_s}
Let $E$ be a closed subspace of $s$ such that $s/E$ has a dominating norm. Then for every dominating Hilbert norm $||\cdot||_E$ on $E$ there is a dominating
Hilbert norm $||\cdot||$ on $s$ such that the restriction of $||\cdot||$ to $E$ is equal to $||\cdot||_E$. 
\end{Cor}
\bproof
The results follows by applying Lemma \ref{lem_extenstion_of_D-Hilbert_N} to the canonical short exact sequence
\[0\longrightarrow E\hookrightarrow{}s\longrightarrow s/E\longrightarrow 0.\]
\bqed

Corollary \ref{cor_extension_hilbertDN_s} together with some results of Vogt (\cite[Cor. 7.6]{V3}, \cite[Prop. 3.3(a)]{V4}) allows to give a stronger version of the well-kown theorem of Vogt characterizing closed 
subspaces of $s$ in terms of nuclearity and the property (DN) (see \cite[Satz 1.7]{V5} or \cite[Prop. 31.5]{MeV}).

%\begin{Th}[Vogt 1977]\label{th_vogt_1977} 
%For a Fr\'echet space $X$ the following assertions are equivalent:
%\begin{itemize}
% \item[\upshape(1)] $X$ is isomorphic to some closed subspace of $s$;
% \item[\upshape(2)] there is a topological embedding $X\hookrightarrow s$;
% \item[\upshape(3)] $X$ is nuclear and has the property $\mathrm{(DN)}$. 
%\end{itemize}
%\end{Th}

\begin{Th}\label{th_top_embedding_DN}
For a Fr\'echet space $X$ and a seminorm $||\cdot||\colon X\to[0,\infty)$ the following assertions are equivalent:
\begin{enumerate}
 \item[\upshape(1)] there is a topological embedding $V\colon X\hookrightarrow s$ such that $||Vx||_{\ell_2}=||x||$ for all $x\in X$;
 \item[\upshape(2)] $X$ is nuclear and $||\cdot||$ is a dominating Hilbert norm on $X$.
\end{enumerate}
\end{Th}

\bproof
($\Rightarrow$) By assumption, $X$ is isomorphic to a closed subspace of $s$, so it is nuclear. Moreover, since $||\cdot||_{\ell_2}$ is a dominating Hilbert norm on $s$, 
it is easily seen that $||\cdot||$ is a dominating Hilbert norm on $X$.

($\Leftarrow$) 
Since $X$ is a nuclear Fr\'echet space with the property (DN), by Theorem \cite[Prop. 31.5]{MeV}, it is isomorphic to some closed subspace $E$ of $s$.
Without lost of generality we may assume that $s/E$ has (DN). Indeed, by \cite[Prop. 3.3(a)]{V4}, for an arbitrary closed subspace $E_0$ of $s$ isomorphic to $X$, 
there is a closed subspace $F_0$ of $s$ and an exact sequence
\[0\longrightarrow E_0\stackrel{\iota_0}\longrightarrow s\stackrel{q_0}\longrightarrow F_0\longrightarrow 0\]
with continuous linear maps.
Then $q_0$ induces an isomorphism $\overline{q_0}\colon s/\ker q_0\to F_0$ (see \cite[Prop. 22.11]{MeV}), and thus 
\[s/\iota_0(E_0)=s/\ker q_0\cong F_0.\]
Hence $E:=\iota_0(E_0)$ is a closed subspace of $s$ such that $s/E$ has (DN), and by the open mapping theorem $E\cong E_0$. 

Let $T\colon X\to E$ be an isomorphism and set $||x||_E:=||T^{-1}x||$ for $x\in E$. 
An easy calculation shows that $||\cdot||_E$ is a dominating Hilbert norm on $E$.  
%Take $q\in\N_0$, $x\in E$ and let $y:=T^{-1}x$. Then, since $||\cdot||$ is a dominating norm on $X$ and $T$ is an isomorphism, 
%we can find $r_1,r_2,r_3\in\N_0$ and constants $C_1,C_2,C_3>0$ such that
%\begin{align*}
%|x|_q^2&=|Ty|^2_q\leq C_1|y|^2_{X,r_1}\leq C_2||y||\;|y|_{X,r_2}=C_2||T^{-1}x||\;|T^{-1}x|_{X,r_2}\\
%&\leq C_3||x||_E|x|_{r_3},
%\end{align*}
%and therefore $||\cdot||_E$ is a dominating norm on $E$. 
Hence, by Corollary \ref{cor_extension_hilbertDN_s}, there is a dominating Hilbert norm $||\cdot||_s$ on $s$ such that ${||\cdot||_s}_{\mid{E}}=||\cdot||_E$.
Moreover, by \cite[Cor. 7.6]{V3}, there is an automorphism $U$ of $s$ such that $||Ux||_{\ell_2}=||x||_s$ for $x\in s$. 
Now, for $V:=U\circ T$ and $x\in X$, we get
\[||Vx||_{\ell_2}=||U(Tx)||_{\ell_2}=||Tx||_s=||Tx||_E=||x||,\]
which proves the theorem.
\bqed

Finally, we are ready to proof our main result.

\proofmain
(3)$\Rightarrow$(1): We should only note that, by \cite[Prop. 4.2]{Cias2}, $\lambda^2(A)\cong\lambda^\infty(|f_j|_q)$ as a Fr\'echet ${}^*$-algebras and $\lambda^\infty(|f_j|_q)$ is isomorphic 
as a Fr\'echet ${}^*$-algebra to the closed ${}^*$-subalgebra of $\ldss$ generated by the sequence $(f_j\otimes f_j)_{j\in\N}$.

(1)$\Rightarrow$(2): Let $E$ be a closed commutative ${}^*$-subalgebra of $\ldss$ isomorphic to the Fr\'echet ${}^*$-algebra $\lambda^2(A)$. By \cite[Th. 4.8]{Cias} and the nuclearity of $E$, there is a
sequence $(Q_j)_{j\in\N}$ of nonzero pairwise orthogonal self-adjoint projections belonging to $\cL(s',s)$ such that
$\lambda^2(A)\cong E\cong\lambda^2(||Q_j||_q)$ as Fr\'echet ${}^*$-algebras. Therefore, according to \cite[Prop. 4.2]{Cias2}, there is a bijection $\sigma\colon\N\to\N$ such that for $P_j:=Q_{\sigma(j)}$ 
($j\in\N$), we have $A\sim(||P_j||_q)_{j\in\N,q\in\N_0}$.

(2)$\Rightarrow$(4): By \cite[Prop. 4.7]{Cias} and \cite[Lemma 27.25]{MeV}, $\lambda^2(||P_j||_q)$ is isomorphic to the closed ${}^*$-subalgebra 
of $\ldss$ generated by the sequence $(P_j)_{j\in\N}$. 
Therefore, the space $\lambda^2(A)$ (which is isomorphic as a Fr\'echet ${}^*$-algebra to $\lambda^2(||P_j||_q)$) is nuclear, and consequently, the matrix $A$ satisfies condition (i).
Moreover, it is easy to check that the matrix $(||P_j||_q)_{j\in\N,q\in\N_0}$ fulfils conditions (ii) and (iii), 
and thus (ii) and (iii) are also satisfied by the matrix $A$ (apply \cite[Prop. 4.2]{Cias2}).

(4)$\Rightarrow$(5): By conditions (ii) and (iii), we get
\[a_{j,r}\leq Ca_{j,p}a_{j,r}\leq Ca_{j,\max\{p,r\}}^2\]
and
\[a_{j,q}^2\leq C_qa_{j,r}\leq CC_qa_{j,p}a_{j,r}\]
with appropriate quantifiers and constants.

(5)$\Leftrightarrow$(6): It is well known that condtions (i) and (ii) are equivalent to nuclearity and the property (DN), respectively.
By \cite[Prop. 3.1, Th. 5.9]{Pir}, $\lambda^\infty(A)$ is a biprojective K\"othe algebra if and only if $A$ satifies condition (iii). 

(5)$\Leftrightarrow$(7): By \cite[Prop. 31.5]{MeV}, conditions (i) and (ii) are equivalent to the fact that $\lambda^2(A)$ is isomorphic to a closed subspace of $s$ and, by \cite[Prop. 3.1]{Pir}, 
$\lambda^\infty(A)$ is an algebra if and only if $A\prec A^2$. Moreover, by \cite[Lemmas 7.7 and 7.9]{Pir}, conditions (i) and (iii) imply that the square root of each nonnegative element of 
$\lambda^\infty(A)$ belongs to $\lambda^\infty(A)$. On the other hand, if $\lambda^2(A)$ is closed under taking square root and nuclear (note that closed subspaces of $s$ are nuclear)
then, by \cite[Prop. 28.16]{MeV}, $\lambda^2(A)\subset\lambda^2(A^2)$.
Consequently, by the closed graph theorem, the identity map from  $\lambda^2(A)$ to $\lambda^2(A^2)$ is continuous, and hence $A^2\prec A$.

(5)$\Rightarrow$(8): Nuclearity is guaranteed by condition (i). For an index $p$ from condition (ii) take, according to (i), an index $r$ and a constant $C_0>0$ 
such that $\sum_{j=1}^\infty\frac{a_{j,p}}{a_{j,r}}<C_0$.
By (ii) and (iii) (more precisely by $A\prec A^2$), there are $q_1,q_2,q_3\in\N_0$ and constants $C_1,C_2,C_3>0$ such that
\[a_{j,r}\leq C_1 a_{j,q_1}^2\leq C_2a_{j,p}a_{j,q_2}\leq C_3a_{j,p} a_{j,q_3}^2.\]
Then for $\xi\in\lambda^2(A)$ we obtain
\[||\xi||_{\ell_2}\leq C_3\sum_{j=1}^\infty|\xi_j|^2\frac{a_{j,p} a_{j,q_3}^2}{a_{j,r}}\leq C_3 \sum_{j=1}^\infty\frac{a_{j,p}}{a_{j,r}}\cdot\sup_{j\in\N}|\xi_j|^2a_{j,q_3}^2
\leq C_0C_3|\xi|_{\lambda^2(A),q_3}^2<\infty,\]
and thus $||\cdot||_{\ell_2}$ is a norm on $\lambda^2(A)$.

%From condition (iv), $a_{j,p}>0$ for all $j\in\N$, so we may set
%$A':=(\frac{a_{j,q}}{a_{j,p}})_{j\in\N,q\in\N_0}$. Now, it is easily seen that $\lambda^\infty(A)=\lambda^\infty(A')$ as Fr\'echet ${}^*$-algebras 
%and $||\cdot||_{\ell_2}$ is a dominating norm on $\lambda^\infty(A')$, so is on $\lambda^\infty(A)$.

Moreover, the Cauchy-Schwartz inequality and condition (iii) (more precisely relation $A^2\prec A$) imply that for all $q\in\N_0$ there are $r\in\N_0$ and $C>0$ such that
\begin{align*}
|\xi|_{\lambda^2(A),q}^2&:=\sum_{j\in\N}|\xi_j|^2a_{j,q}^2\leq C\bigg(\sum_{j\in\N}|\xi_j|^2\bigg)^{1/2}\cdot\bigg(\sum_{j\in\N}|\xi_j|^2 a_{j,q}^4\bigg)^{1/2}
\leq C||\xi||_{\ell_2}\left(\sum_{j\in\N}|\xi_j|^2 a_{j,r}^2\right)^{1/2}\\
&= C||\xi||_{\ell_2}|\xi|_{\lambda^2(A),r}, 
\end{align*}
which shows that $||\cdot||_{\ell_2}$ is a dominating norm on $\lambda^2(A)$.

%Moreover, condition (ii) and the Schwartz inequality imply that for all $q\in\N_0$ there are $r\in\N_0$ and $C_q>0$ such that
%\begin{align*}
%|\xi|_{\lambda^2(A),q}^2&:=\sum_{j\in\N}|\xi_j|^2a_{j,q}^2\leq C_q\sum_{j\in\N}|\xi_j|^2a_{j,r}
%\leq C_q\bigg(\sum_{j=1}^\infty|\xi_j|^2\bigg)^{1/2}\bigg(\sum_{j=1}^\infty|\xi_j|^{2}a_{j,r}^{2}\bigg)^{1/2}\\
%&=C_q||\xi||_{\ell_2}|\xi|_{\lambda^{2}(A),r}.\\
%\end{align*}
%This shows that $||\cdot||_{\ell_2}$ is a dominating norm on $\lambda^2(A)$. 

(8)$\Rightarrow$(3): 
By Theorem \ref{th_top_embedding_DN}, there is a topological embedding $V\colon\lambda^2(A)\hookrightarrow s$ such that 
$||V\xi||_{\ell_2}=||\xi||_{\ell_2}$ for all $\xi\in \lambda^2(A)$.
%and thus we may define $\widetilde{V}\colon\ell_2\to\ell_2$ by $\widetilde{V}x:=\lim_{j\to\infty}Vx_j$, where $(x_j)_{j\in\N}$ is an arbitrary sequence in $\lambda^2(A)$ such that $\lim_{j\to\infty}x_j=x$.
%Clearly, $||\widetilde{V}x||_{\ell_2}=||x||_{\ell_2}$ for $x\in\ell_2$ and $\widetilde{V}$ has dense range, hence $\widetilde{V}$ is a unitary map.
Set $f_j:=Ve_j$ for $j\in\N$. Then $(f_j)_{j\in\N}\subset s$ is an orthonormal sequence in $\ell_2$ and a Schauder basis of $\im V$. 
Therefore, $\Phi\colon\im V\to\lambda^2(|f_j|_q)$ defined by $\Phi f_j:=e_j$ for $j\in\N$ is a Fr\'echet space isomorphism (see \cite[Lemma 27.25]{MeV}),
and so is $\Phi\circ V\colon\lambda^2(A)\to\lambda^2(|f_j|_q)$. 
But $(\Phi\circ V)e_j=e_j$ for $j\in\N$, hence $\lambda^2(A)=\lambda^2(|f_j|_q)$ as Fr\'echet ${}^*$-algebras, whence $A\sim(|f_j|_q)_{j\in\N,q\in\N_0}$.
\bqed

%Since $\lambda^\infty(A)$ has the property (DN) it is isomorphic as a Fr\'echet space to some closed subspace $F$ of $s$. 

%By nuclearity and \cite[Lemma 31.4]{MeV}, there is a closed subspace $F$ of $s$ and an exact sequence
%\[0\longrightarrow s\stackrel{j}\longrightarrow F\stackrel{q}\longrightarrow\lambda^\infty(A)\longrightarrow 0\]
%with continuous linear maps between Fr\'echet spaces.
%Since $\lambda^\infty(A)$ has the property (DN), the above sequence splits (for details, see \cite[Prop. 31.5]{MeV} and its proof), i.e. $q$ has a continuous linear right inverse $u\colon\lambda^\infty(A)\to F$.
%In consequence, $\lambda^\infty(A)$ is isomorphic to a closed subspace of $F$.

%\begin{problem}
%Which K\"othe matrices $A$ admit a topological embedding $u\colon\lambda^2(A)\to s$ such that:
%\begin{itemize}
% \item $u$ extends to unitary operator $\widetilde{u}:$
%\end{itemize}  
%\end{problem}

\section{Quasi-equivalence Conjecture}

In the last section we shall show how commutative subalgebras of $\ldss$ are connected with the so-called Quasi-equivalence Conjecture stating that all bases in a nuclear Fr\'echet space are quasi-equivalent
(see \cite[\S 8.3]{Mit}, \cite{Zobin} and references therein).
Here we only consider the case (still unsolved) of closed subspaces of $s$, i.e. nuclear Fr\'echet spaces with the property (DN).
Let us recall that two bases $(f_j)_{j\in\N}$ and $(g_j)_{j\in\N}$ of a Fr\'echet space $X$ are called \emph{quasi-eqivalent} if  
there is a bijection $\sigma\colon\N\to\N$ and a sequence $(\lambda_j)_{j\in\N}$ of nonzero scalars such that 
the operator $T\colon X\to X$ defined by $Tf_j=\lambda_jg_{\sigma(j)}$ is a Fr\'echet space isomorphism.

\begin{Th}\label{th_quasi-eq}
Let $X$ be an infinite-dimensional closed subspace of $s$.
Then the following assertions are equivalent.
\begin{enumerate}
\item[\upshape(i)] For each two closed commutative ${}^*$-subalgebras $E,F$ of $\ldss$ we have: if $E\cong F\cong X$ as Fr\'echet spaces then $E\cong F$ as Fr\'echet ${}^*$-algebras. 
\item[\upshape(ii)] For each two basic sequences $(f_j)_{k\in\N}$ and $(g_j)_{k\in\N}$ of $s$ which are at the same time orthonormal sequences in $\ell_2$ we have: 
if $\lambda^2(|f_j|_q)\cong\lambda^2(|g_j|_q)\cong X$ as Fr\'echet spaces then there is a bijection $\sigma\colon\N\to\N$ such that $(|f_j|_q)_{j\in\N,q\in\N_0}\sim (|g_{\sigma(j)}|_q)_{j\in\N,q\in\N_0}$.
 %\begin{itemize}
 %\item[($\alpha$)] $\forall q\in\N_0\;\exists r\in\N_0\;\exists C>0\;\forall j\in\N\quad|f_j|_q\leq C|g_{\sigma(j)}|_r$;
 %\item[($\beta$)] $\forall r'\in\N_0\;\exists q'\in\N_0\;\exists C'>0\;\forall j\in\N\quad |g_{\sigma(j)}|_{r'}\leq C'|f_j|_{q'}$.
 %\end{itemize}
\item[\upshape(iii)] For each two K\"othe matrices $A$ and $B$ such that $A\sim A^2$ and $B\sim B^2$ we have: 
if $\lambda^2(A)\cong\lambda^2(B)\cong X$ as Fr\'echet spaces then there is a bijection $\sigma\colon\N\to\N$ such that $A\sim B_{\sigma}$. 
\item[\upshape(iv)] In $X$ all bases are quasi-equivalent. 
\end{enumerate}  
\end{Th}

\bproof
(i)$\Leftrightarrow$(ii): If $E$ and $F$ are closed commutative ${}^*$-subalgebras of $\ldss$ then, by Theorem \ref{th_main}, there are basic sequences $(f_j)_{j\in\N}$ and $(g_j)_{j\in\N}$ of $s$
which are at the same time orthonormal in $\ell_2$ such that $E\cong\lambda^2(|f_j|_q)$ and $F\cong\lambda^2(|g_j|_q)$ as Fr\'echet ${}^*$-algebras. Moreover, for each such a sequence
$(f_j)_{j\in\N}$, $\lambda^2(|f_j|_q)$ is isomorphic as a Fr\'echet ${}^*$-algebra to a closed commutative ${}^*$-subalgebra of $\ldss$.
Therefore, assuming $\lambda^2(|f_j|_q)\cong\lambda^2(|g_j|_q)\cong X$ as Fr\'echet spaces, it follows from (i) and \cite[Prop. 4.2]{Cias2} that 
$(|f_j|_q)_{j\in\N,q\in\N_0}\sim (|g_{\sigma(j)}|_q)_{j\in\N,q\in\N_0}$. 
Conversely, if we assume that $E\cong F\cong X$ as Fr\'echet spaces then, by (ii) and again by \cite[Prop. 4.2]{Cias2}, we get that $E\cong F$ as Fr\'echet ${}^*$-algebras.

(i)$\Leftrightarrow$(iii): This is an immediate consequence of Theorem \ref{th_main} and \cite[Prop. 4.2]{Cias2}.

(iii)$\Rightarrow$(iv): Let $(x_j)_{j\in\N}$, $(y_j)_{j\in\N}$ be two bases in $X$. Then 
\[\lambda^2(A)\cong\lambda^2(B)\cong X\]
as Fr\'echet spaces, where 
\[A:=\bigg(\frac{|x_j|_q}{||x_j||_{\ell_2}}\bigg)_{j\in\N,q\in\N_0}\quad \text{and} \quad B:=\bigg(\frac{|y_j|_q}{||y_j||_{\ell_2}}\bigg)_{j\in\N,q\in\N_0}.\]
Then for all $q\in\N_0$ we have
\[\bigg(\frac{|x_j|_q}{||x_j||_{\ell_2}}\bigg)^2\leq\frac{||x_j||_{\ell_2}|x_j|_{2q}}{||x_j||_{\ell_2}^2}=\frac{|x_j|_{2q}}{||x_j||_{\ell_2}}\]
so $A^2\prec A$. Moreover, since
\[\frac{|x_j|_q}{||x_j||_{\ell_2}}\geq1\]
for all $q\in\N_0$, we have also $A\prec A^2$, and thus $A\sim A^2$. Similarly, we show that $B\sim B^2$. 
Therefore, by (iii), there is a bijection $\sigma\colon\N\to\N$ such that $A\sim B_\sigma$. This means that
 \[\forall q\in\N_0\;\exists r\in\N_0\;\exists C>0\;\forall j\in\N\quad|x_j|_q\leq C\lambda_j|y_{\sigma(j)}|_r\]
and
 \[\forall r'\in\N_0\;\exists q'\in\N_0\;\exists C'>0\;\forall j\in\N\quad \lambda_j|y_{\sigma(j)}|_{r'}\leq C'|x_j|_{q'},\]
where $\lambda_j:=\frac{||x_j||_{\ell_2}}{||y_{\sigma(j)}||_{\ell_2}}$. Hence $x_j\mapsto \lambda_jy_{\sigma(j)}$, $j\in\N$, defines an isomorphism of $X$, i.e. the bases $(x_j)_{j\in\N}$ and $(y_j)_{j\in\N}$ 
are quasi-equivalent. 

(iv)$\Rightarrow$(ii): Assume that in $X$ all bases are quasi-equivalent. Let $(f_j)_{j\in\N}$ and $(g_j)_{j\in\N}$ be basic sequences of $s$ which are at the same time orthonormal sequences of $\ell_2$ 
and such that $\lambda^2(|f_j|_q)\cong\lambda^2(|g_j|_q)\cong X$ as Fr\'echet spaces.  Let $F$ and $G$ be closed linear span in $s$ of $(f_j)_{j\in\N}$ and $(g_j)_{j\in\N}$, respectively.
Then $F\cong G\cong X$ as Fr\'echet spaces; let $T\colon G\to F$ be an isomorphism. 
Clearly, $(F_j)_{j\in\N}$, $F_j:=Tg_j$, is a basis in $F$ which is, by assumption quasi-equivalent to $(f_j)_{j\in\N}$. Consequently,
there is a sequence $(\lambda_j)_{j\in\N}$ of nonzero scalars and a bijection $\sigma\colon\N\to\N$ such that $S\colon F\to F$ defined by sending $F_{\sigma(j)}$ to $\lambda_j^{-1}f_j$ is 
an isomorphism of Fr\'echet spaces.

Now, let $V:=ST\colon G\to F$ and define $||\cdot||_V\colon F\to[0,\infty)$, by $||V\xi||_V:=||\xi||_{\ell_2}$. Since $||\cdot||_{\ell_2}$ is a dominating Hilbert norm on $G$, $||\cdot||_V$
is a dominating Hilbert norm on $F$ and, consequently, $(\lambda_j^{-1}f_j)_{j\in\N}=(Vg_{\sigma(j)})_{j\in\N}$ is an orthonormal sequence in $F$ with respect to the Hilbert norm $||\cdot||_V$. 
In particular, $||f_j||_V=|\lambda_j|$ for $j\in\N $ and thus, $||\cdot||_V$ being a dominating norm on $F$, 
\begin{equation}\label{eq_1lambda_j}
\forall q\in\N_0\;\exists r\in\N_0,C>0\;\forall j\quad |f_j|_q\leq C||f_j||_V|f_j|_r=C|\lambda_j|\;|f_j|_r. 
\end{equation}
Moreover, since $||\cdot||_V$ is a continuous norm on $F$,
\[\exists r_0\in\N_0,C_0>0\;\forall j\quad |\lambda_j|=||f_j||_V\leq C_0|f_j|_{r_0}.\]
Now, since
\[\forall q\in\N_0\;\exists r\in\N_0,C>0\;\forall j\in\N\quad |f_j|_q^2\leq C|f_j|_r,\]
it follows that
\begin{equation}\label{eq_2lambda_j}
\forall q\in\N_0\;\exists r\in\N_0,C>0\;\forall j\quad |\lambda_j|\;|f_j|_q\leq C_0|f_j|_{r_0}|f_j|_q\leq C_0|f_j|_{\max\{r_0,q\}}^2\leq C|f_j|_r.
\end{equation}

Finally, let $W\colon F\to F$ be defined by $f_j\mapsto\lambda_jf_j$ for $j\in\N$. Then, by (\ref{eq_1lambda_j}) and (\ref{eq_2lambda_j}), $W$ is an automorphism of the Fr\'echet space $F$,
and thus $WV$ is an isomorphism of Fr\'echet spaces which sends $g_{\sigma(j)}$ to $f_j$. This, clearly, implies that $(|f_j|_q)_{j\in\N,q\in\N_0}\sim (|g_{\sigma(j)}|_q)_{j\in\N,q\in\N_0}$.  
\bqed

For a monotonically increasing sequence $\alpha=(\alpha_j)_{j\in\N}$ in $[0,\infty)$ such that $\lim_{j\to\infty}\alpha_j=\infty$ we define the \emph{power series space of infinite type}
\[\Lambda_\infty(\alpha):=\{(\xi_j)_{j\in\N}\subset\C^\N\colon \sum_{j=1}^\infty|\xi_j|^2e^{2q\alpha_j}<\infty\quad\text{for all $q\in\N_0$}\}.\]
It appears that the space $\Lambda_\infty(\alpha)$ is nuclear if and only if $\sup_{j\in\N}\frac{\log j}{\alpha_j}<\infty$ (see \cite[Prop. 29.6]{MeV}).
As a consequence, of Theorem \ref{th_quasi-eq}, we get another proof of the following well-known result of Mityagin \cite[Th. 12]{Mit}.
\begin{Cor}[Mityagin 1961]
In nuclear power series space of infinite type $\Lambda_\infty(\alpha)$ all bases are quasi-equivalent.
\end{Cor}

\bproof
Let $E, F$ be a closed commutative ${}^*$-subalgebras of $\ldss$  isomorphic as Fr\'echet spaces to $\Lambda_\infty(\alpha)$. 
Then, by \cite[Cor. 6.10]{Cias2}, $E$ and $F$ are isomorphic to $\Lambda_\infty(\alpha)$ as a Fr\'echet ${}^*$-algebra, and thus, by Theorem \ref{th_quasi-eq}, all bases in $\Lambda_\infty(\alpha)$ 
are quasi-equivalent.
\bqed

Since, by Corollary \ref{cor_frechet_spaces_subalg_ldss}, every closed subspace of $s$ with basis is isomorphic as a Fr\'echet space to some closed commutative ${}^*$-subalgebra of $\ldss$, 
we can rewrite Theorem \ref{th_quasi-eq} in the following way.

\begin{Cor}\label{cor_quasi-eq2}
The following assertions are equivalent.
\begin{enumerate}
\item[\upshape(i)] For each two closed commutative ${}^*$-subalgebras $E,F$ of $\ldss$ we have: if $E\cong F$ as Fr\'echet spaces then $E\cong F$ as Fr\'echet ${}^*$-algebras. 
\item[\upshape(ii)] For each two basic sequences $(f_j)_{k\in\N}$ and $(g_j)_{k\in\N}$ of $s$ which are at the same time orthonormal sequences in $\ell_2$ we have: 
if $\lambda^2(|f_j|_q)\cong\lambda^2(|g_j|_q)$ as Fr\'echet spaces then there is a bijection $\sigma\colon\N\to\N$ such that $(|f_j|_q)_{j\in\N,q\in\N_0}\sim (|g_{\sigma(j)}|_q)_{j\in\N,q\in\N_0}$.
 %\begin{itemize}
 %\item[($\alpha$)] $\forall q\in\N_0\;\exists r\in\N_0\;\exists C>0\;\forall j\in\N\quad|f_j|_q\leq C|g_{\sigma(j)}|_r$;
 %\item[($\beta$)] $\forall r'\in\N_0\;\exists q'\in\N_0\;\exists C'>0\;\forall j\in\N\quad |g_{\sigma(j)}|_{r'}\leq C'|f_j|_{q'}$.
 %\end{itemize}
\item[\upshape(iii)] For each two K\"othe matrices $A$ and $B$ such that $A\sim A^2$ and $B\sim B^2$ we have: 
if $\lambda^2(A)\cong\lambda^2(B)$ as Fr\'echet spaces then there is a bijection $\sigma\colon\N\to\N$ such that $A\sim B_{\sigma}$. 
\item[\upshape(iv)] In each closed subspace of $s$ all bases are quasi-equivalent. 
\end{enumerate}  
\end{Cor}

%$\begin{Conj}
%Every closed subspace of $s$ with basis has the quasi-equivalence property.
%\end{Conj}

\begin{Acknow*}
I cordially thank Leonhard Frerick for many valuable suggestions and encouraging me to solve the problems raised in the present paper.     
\end{Acknow*}

\vspace{1cm}
\begin{minipage}{8.5cm}
T. Cia\'s

Faculty of Mathematics and Computer Science

Adam Mickiewicz University in Pozna{\'n}

Umultowska 87

61-614 Pozna{\'n}, POLAND

e-mail: tcias@amu.edu.pl
\end{minipage}\

\end{document}